\documentclass[12pt]{article}
\newcommand{\di}{\displaystyle}
\usepackage{amsmath}
\usepackage{amssymb}
\usepackage{amsthm}

\newtheorem{definition}{Definition}
\newtheorem{theorem}{Theorem}

\newtheorem{coro}{Corollary}
\begin{document}

\title{Rarita-Schwinger Type Operators on Spheres and Real Projective Space}

\author{Junxia Li and John Ryan\\
\emph{\small Department of Mathematics, University of Arkansas, Fayetteville, AR 72701, USA}\\
Carmen J. Vanegas\\
\emph{\small Departamento de Matem\'aticas, Universidad Sim\'on Bol\'ivar, Caracas, Venezuela}}

\date{}

\maketitle
\begin{abstract}
In this paper we deal with Rarita-Schwinger type operators on spheres and real projective space.
First we define the spherical Rarita-Schwinger type operators and construct their fundamental solutions.
Then we  establish that the projection operators appearing in the spherical Rarita-Schwinger type operators and
the spherical Rarita-Schwinger type equations are conformally invariant under the Cayley transformation.
Further, we obtain some basic integral formulas related to the spherical Rarita-Schwinger type operators.
Second, we define the Rarita-Schwinger type operators on the real projective space and
construct their kernels and Cauchy integral formulas.
\end{abstract}

\section{Introduction}

 Rarita-Schwinger operators are generalizations of the Dirac operator and arise in representation theory for the Spin and Pin groups. See \cite{BSSV1, BSSV2,DLRV, Va1, Va2}. We denote a Rarita-Schwinger operator by $R_k$, where $k=0, 1, \cdots, m, \cdots.$ When $k=0$ we have the Dirac operator. The Rarita-Schwinger operators $R_k$ in Euclidean space have been studied in \cite{BSSV1, BSSV2,DLRV, Va1, Va2}. Here we construct similar Rarita-Schwinger operators together with their fundamental solutions and study their representation theory on the sphere and real projective space.

First J. Ryan \cite{R1,R2} in 1997 and P. Van Lancker \cite{Va3} in 1998 studied the Dirac operators on the sphere. Later, H. Liu and J. Ryan \cite{LR} studied the spherical Dirac type operators on the sphere by using Cayley transformations. See also \cite{BR}. Using similar methods to define the Rarita-Schwinger operators in $\mathbb{R}^n$, we can define the spherical Rarita-Schwinger type operator on the sphere based on the spherical Dirac operator. We also use similar arguments as in Euclidean space to establish the conformal invariance for the projection operators and the spherical Rarita-Schwinger type equations under the Cayley transformations. See \cite{DLRV}. Further the fundamental solutions to the spherical Rarita-Schwinger type operators are achieved by applying the Cayley transformation. In turn, Stokes' Theorem, Cauchy's Theorem, Borel-Pompeiu Formula, Cauchy Integral Formula and a Cauchy Transform are proved for the sphere. Furthermore, we show that Stokes' theorem is conformally invariant under Cayley transformation, and with minor modification, is equivalent to the Rarita-Schwinger version of Stokes' Theorem in Euclidean space appearing in \cite{BSSV1, DLRV} and elsewhere.

 By factoring out $\mathbb{S}^n$ by the group $\mathbb{Z}_2=\{\pm 1\}$ we obtain real projective space, $\mathbb{R}P^n$. On this space, we define the Rarita-Schwinger type operators and construct their kernels over two different bundles over $\mathbb{R}P^n$. Further, we obtain some basic integral formulas from Clifford analysis associated with these operators for the two different bundles. This extends results from \cite{KR}.

\section{Preliminaries}

 A Clifford algebra, $Cl_{n+1},$ can be generated from $\mathbb{R}^{n+1}$ by considering the
relationship $$\underline{x}^{2}=-\|\underline{x}\|^{2}$$ for each
$\underline{x}\in \mathbb{R}^{n+1}$.  We have $\mathbb{R}^{n+1}\subseteq Cl_{n+1}$. If $e_1,\ldots, e_{n+1}$ is an orthonormal basis for $\mathbb{R}^{n+1}$, then $\underline{x}^{2}=-\|\underline{x}\|^{2}$ tells us that $e_i e_j + e_j e_i= -2\delta_{ij}.$ Let $A=\{j_1, \cdots, j_r\}\subset \{1, 2, \cdots, n+1\}$ and $1\leq j_1< j_2 < \cdots < j_r \leq n+1$. An arbitrary element of the basis of the Clifford algebra can be written as $e_A=e_{j_1}\cdots e_{j_r}.$ Hence for any element $a\in Cl_{n+1}$, we have $a=\sum_Aa_Ae_A,$ where $a_A\in \mathbb{R}.$

The reversion is given by
 $$\tilde{a}=\sum_A(-1)^{|A|(|A|-1)/2}a_Ae_A,$$ where $|A|$ is the cardinality of $A$. In particular, $\widetilde{e_{j_1}\cdots e_{j_r}}=e_{j_r}\cdots e_{j_1}.$ Also $\widetilde{ab}=\tilde{b}\tilde{a}$ for $a, b \in Cl_{n+1}.$
The Clifford conjugation is defined by
 $$\bar{a}=\sum_A(-1)^{|A|(|A|+1)/2}a_Ae_A$$
and satisfies $\overline{e_{j_1}\cdots e_{j_r}}=(-1)^r e_{j_r}\cdots e_{j_1}$ and $\overline{ab}= \bar{b}\bar{a}$ for $a, b \in Cl_{n+1}.$

For each $a=a_0+a_1e_1+\cdots +a_{1\cdots n+1}e_1\cdots e_{n+1}\in Cl_{n+1}$ the scalar part of $\bar{a}a$ gives
the square of the norm of $a,$ namely $a_0^2+a_1^2+\cdots +a_{1\cdots n+1}^2$\,.
For more on Clifford algebras and their properties, see \cite{P}.

The Pin and Spin groups play an important role in Clifford analysis.  The Pin group can be defined as
 $$Pin(n+1): =\{a\in Cl_{n+1} : a=y_1 \ldots y_p:
{y_1,\ldots , y_p}\in \mathbb{S}^{n}, p\in \mathbb{N}\}$$
and it is clearly a group under multiplication in $Cl_{n+1}$.
\par Now suppose that $y\in \mathbb{S}^{n}\subseteq \mathbb{R}^{n+1}$. Look at $yxy=yx^{\parallel _y}y+yx^{\perp_y}y=-x^{\parallel _y}+x^{\perp_y}$ where $x^{\parallel _y}$ is the projection of $x$ onto $y$
and $x^{\perp_y}$ is perpendicular to $y$. So $yxy$ gives a reflection of $x$ in the $y$ direction. By the Cartan$-$Dieudonn\'{e} Theorem each $O \in O(n+1)$  is the composition of a finite number of reflections. If  $a=y_1\ldots y_p\in Pin(n+1)$, then $\tilde{a}:=y_p\ldots y_1$ and $ax\tilde{a}=O_a(x)$ for some $O_a\in O(n+1).$ Choosing $y_1, \ldots, y_p$ arbitrarily in $\mathbb{S}^{n}$,  we see that the group homomorphism

$$\theta: Pin(n+1)\longrightarrow O(n+1): a\longmapsto O_a,$$
 with $a=y_1\ldots y_p$
and $O_a(x)=ax\tilde{a},$ is surjective. Further $-ax(-\tilde{a})=ax\tilde{a}$, so $1, -1\in ker(\theta)$. In fact $ker(\theta)=\{\pm 1\}.$ See \cite{P}.
The Spin group is defined as
$$
Spin(n+1):=\{a\in Pin(n+1): a=y_1\ldots y_ p \mbox{ and } p \mbox{ even}\}
$$
and it is a subgroup of $Pin(n+1)$. There is a group homomorphism
 $$\theta: Spin(n+1)\longrightarrow SO(n+1)$$ which is surjective with kernel $\{1, -1\}$. See \cite{P} for details.

The Dirac Operator in $\mathbb{R}^{n}$ is defined to be $$D :=\sum_{j=1}^{n} e_j \frac{\partial}{\partial x_j}.$$ Note $D^2=-\Delta_{n},$ where $\Delta_{n}$ is the Laplacian in $\mathbb{R}^{n}$.

If $p_k$ is a homogeneous polynomial with degree $k$ such that $Dp_k=0,$ we call such a polynomial a left monogenic polynomial homogeneous of degree $k$.

Let $\mathcal{H}_{k}$ be the space of $Cl_{n+1}-$ valued harmonic polynomials homogeneous of degree $k$
and $\mathcal{M}_k$ the space of $Cl_{n+1}-$ valued monogenic polynomials homogeneous of degree $k$.
Note if $h_k\in$ $\mathcal{H}_k,$ then $Dh_k\in$ $\mathcal{M}_{k-1}$. But $Dup_{k-1}(u)=(-n-2k+2)p_{k-1}(u),$
so
$$
\mathcal{H}_k=\mathcal{M}_k\bigoplus u\mathcal{M}_{k-1}, h_k=p_k+up_{k-1}.
$$
This is the so-called Almansi-Fischer decomposition of $\mathcal{H}_k$. See \cite{BDS, R}.

Note that if $ p(u)\in \mathcal{M}_k $ then it trivially extends to $P(v) = p(u + u_{n+1}e_{n+1})$
with $u_{n+1}\in \mathbb{R}$ and $P(v)= p(u)$ for all $u_{n+1}\in \mathbb{R}$. Consequently $D_{n+1}P(v)= 0$
where $D_{n+1}= \displaystyle\sum_{j=1}^{n+1} {e_j \frac{\partial}{\partial u_j}}$.

If $ p(u)\in \mathcal{M}_k $ then for any boundary of a piecewise smooth bounded domain
$U \subseteq \mathbb{R}^{n} $ by Cauchy's Theorem
\begin{equation}\label{CTu}
\int_{\partial U}{n(u)p(u)d\sigma(u)} = 0.
\end{equation}

Suppose now $a\in \mbox{Pin}(n+1)$ and $u= aw\tilde{a}$ then although $u \in \mathbb{R}^{n}$ in
general $w$ belongs to the hyperplane $a^{-1}\mathbb{R}^{n}\tilde{a}^{-1}$ in $\mathbb{R}^{n+1}$.

By applying a change of variable, up to a sign the integral (\ref{CTu}) becomes
\begin{equation}\label{CVCTu}
\int_{a^{-1}\partial U\tilde{a}^{-1}}{ an(w)\tilde{a}P(aw\tilde{a})d\sigma(w)} = 0.
\end{equation}
As $\partial U$ is arbitrary then on applying Stokes' Theorem to (\ref{CVCTu}) we see
that
\begin{equation}\label{Diraca}
D_a \tilde{a}P(aw\tilde{a}) = 0,~~ \mbox{where} ~~ D_a : = D_{n+1}\bigl|_{a^{-1}\mathbb{R}^{n}\tilde{a}^{-1}}.
\end{equation}

Suppose $U$ is a domain in $\mathbb{R}^{n}$. Consider a function of two variables
$$
f: U\times \mathbb{R}^{n}\longrightarrow Cl_{n+1}
$$
such that for each $x\in U, f(x,u)$ is a left monogenic polynomial homogeneous of degree $k$ in $u$.
Let $P_k$ be the left projection map
$$P_k:  \mathcal{H}_k\rightarrow \mathcal{M}_k,$$ then
$R_kf(x,u)$ is defined to be $P_kD_xf(x,u)$. The left Rarita-Schwinger equation is defined to be
$$
R_k f(x,u)=0.
$$
We also have a right projection $P_{k,r}: \mathcal{H}_k \rightarrow \overline{\mathcal{M}_k},$ and a right Rarita-Schwinger equation $f(x,u)D_xP_{k,r}=f(x,u)R_k=0,$
where $\overline{\mathcal{M}_k}$ stands for the space of right monogenic polynomials homogeneous of degree $k$. See \cite{DLRV}.

\section{Rarita-Schwinger type operators on spheres}
Let $\mathbb{R}^n$ be the span of $e_1,\cdots,e_n.$ Consider the Cayley transformation
$ C : \mathbb{R}^n \to \mathbb{S}^n $, where $\mathbb{S}^n$ is the unit sphere in $\mathbb{R}^{n+1}$,
defined by $C(x)= (e_{n+1}x +1)(x + e_{n+1})^{-1} $, where $ x= x_1e_1 + \cdots + x_ne_n \in \mathbb{R}^n $,
and $e_{n+1}$ is a unit vector in $\mathbb{R}^{n+1}$ which is orthogonal to $\mathbb{R}^n $.
Now $C(\mathbb{R}^n) = \mathbb{S}^n \setminus\{e_{n+1}\}$. Suppose $ x_s \in \mathbb{S}^n $ and
$x_s = x_{s_1}e_1 + \cdots + x_{s_n}e_n + x_{s_{n+1}}e_{n+1}$, then we have
$x=C^{-1}(x_s) = (-e_{n+1}x_s +1)(x_s - e_{n+1})^{-1}$.

The Dirac operator over the $n$-sphere $\mathbb{S}^n$ has the form $D_{s} = w(\Lambda + \frac{n}{2})$, where
$w \in \mathbb{S}^n $ and $ \Lambda =
\displaystyle\sum_{i<j, i=1}^{n}{e_ie_j(w_i\frac{\partial}{\partial w_j} - w_j\frac{\partial}{\partial w_i})}$,
see for instance \cite{CM,LR,Va3}.

Let $U$ be a domain in $\mathbb{R}^n$. Consider a function $f_{\star}: U \times \mathbb{R}^{n} \to Cl_{n+1}$ such that for each $x\in U$,
$f_{\star}(x, u)$ is a left monogenic polynomial homogeneous of degree $k$ in $u$.
This function reduces to $f(x_s, u)$ on $C(U) \times \mathbb{R}^{n} $ and $f(x_s, u)$ takes its values in $Cl_{n+1}$
where $x=C^{-1}(x_s)$ and $x_s\in C(U)\subset \mathbb{S}^n.$ Further
$f(x_s, u)$ is a left monogenic polynomial homogeneous of degree $k$ in $u$.

Since $\Delta_u D_{s, x_s} = D_{s, x_s} \Delta_u$, then
$D_{s, x_s}f(x_s, u)$ is harmonic in $u$. Hence
by the Almansi-Fischer decomposition:
$$
D_{s, x_s}f(x_s, u) = f_{1,k}(x_s, u) + u f_{2,k-1}(x_s, u),
$$
where $f_{1,k}(x_s, u)$
is a left monogenic polynomial homogeneous of degree $k$ in $u$ and
$f_{2,k-1}(x_s, u)$ is a left monogenic polynomial homogeneous of degree $k-1$ in $u$.

We can also consider a function $g_{\star}: U \times \mathbb{R}^{n} \to Cl_{n+1}$ such that for each $x\in U$,
$g_{\star}(x, u)$ is a right monogenic polynomial homogeneous of degree $k$ in $u$.
This function also reduces to a right monogenic polynomial homogeneous $g(x_s, u)$ on $C(U) \times \mathbb{R}^{n} $.

Let $P_k$ be the left projection map
$ P_k : {\cal H}_k = {\mathcal M}_k \oplus u{\mathcal M}_{k-1} \to {\mathcal M}_k$, then the
$n$-spherical left Rarita-Schwinger type operator $ R_k^{S}$  is defined to be
$$
R_k^{S}f(x_s, u) = P_kD_{s,x_s}f(x_s, u).
$$
On the other hand, the $n$-spherical right Rarita-Schwinger type operator $R_{k,r}^{S}$  is defined to be
$$
g(x_s, u)R_{k,r}^{S} = g(x_s, u)D_{s,x_s}P_{k,r},
$$
where $P_{k,r}$ is the right projection
$ P_{k,r} : {\cal H}_k  \to \bar{\mathcal M}_k$.
Consequently, the left and the right $n$-spherical Rarita-Schwinger type equations are defined to be
$$
 R_k^{S}f(x_s, u) = 0 ~~ \mbox{and} ~~  g(x_s, u)R_k^{S} = 0 ~~ \mbox{respectively}.
$$

\section{Conformal invariance of $P_k$ under the Cayley transformation and its inverse}

Consider the Cayley transformation $C(x)= (e_{n+1}x +1)(x + e_{n+1})^{-1}=e_{n+1}(x-e_{n+1})(x + e_{n+1})^{-1}=
e_{n+1}(x+e_{n+1}-2e_{n+1})(x + e_{n+1})^{-1}=e_{n+1}+2(x + e_{n+1})^{-1}.$ This last term is the Iwasawa decomposition
for the Cayley transformation, $C$. Further, $C^{-1}(x_s)= (-e_{n+1}x_s +1)(x_s- e_{n+1})^{-1}=
-e_{n+1}(x_s+e_{n+1})(x - e_{n+1})^{-1}=-e_{n+1}(x_s-e_{n+1}+2e_{n+1})(x_s- e_{n+1})^{-1}=
-e_{n+1}+2(x_s- e_{n+1})^{-1},$ and this last term is the Iwasawa decomposition for the inverse, $C^{-1},$ of the Cayley transformation.

Now let $f(x_s,u): U_s \times \mathbb{R}^{n}\to Cl_{n+1}$ be a monogenic polynomial
homogeneous of degree $k$ in $u$ for each $x_s \in U_s$, where $U_s$ is a domain in $ \mathbb{S}^n.$

It is shown in \cite{DLRV} that $P_k$ is conformally invariant under a general M\"{o}bius transformation
over $\mathbb{R}^{n}$. This trivially extends to M\"{o}bius transformations on $\mathbb{R}^{n+1}$.
It follows that if we restrict $x_s $ to $ \mathbb{S}^n,$ then $P_k$ is also conformally invariant under
the Cayley transformation $C$ and its inverse $C^{-1},$ with $x\in \mathbb{R}^n$.

It follows that we have:
\begin{theorem}
$$
P_{k,w}J(C,x)
f(C(x),\frac{(x+e_{n+1})w(x+ e_{n+1})}{||x+e_{n+1}||^2})=J(C,x)P_{k,u}f(x_s,u),
$$
where
$u=\displaystyle\frac{(x+ e_{n+1})w(x+e_{n+1})}{||x+e_{n+1}||^2}$
and $J(C,x)=\di\frac{x+e_{n+1}}{\|x+e_{n+1}\|^n}$ is the conformal weight for the Cayley transformation.
\end{theorem}
Also for $U$ a domain in $ \mathbb{R}^n,$ and $g(x,u)$ defined on $U \times \mathbb{R}^{n}$
such that for each $x\in U$,
$g$ is monogenic in $u$ and homogeneous of degree $k$ in $u$, we have:
\begin{theorem}\label{t2}
$$
P_{k,w}J(C^{-1},x_s)g(C^{-1}(x_s),
\frac{(x_s - e_{n+1})w(x_s - e_{n+1})}{||x_s - e_{n+1}||^2})=J(C^{-1},x_s)P_{k,u}g(x,u),
$$
where
$u=\displaystyle\frac{(x_s - e_{n+1})w(x_s - e_{n+1})}{||x_s - e_{n+1}||^2}$
and $J(C^{-1},x_s)=\di\frac{x_s - e_{n+1}}{\|x_s - e_{n+1}\|^n}$ is the conformal weight for the inverse Cayley transformation.
\end{theorem}
Note that in the previous  the\-o\-rems $a_1(x):= \di\frac{x+e_{n+1}}{\|x+e_{n+1}\|}$
and $a_2(x_s):=\di\frac{x_s - e_{n+1}}{\|x_s - e_{n+1}\|}$ belong to $\mbox{Pin}(n+1)$.
So $w \in \mathbb{R}^{n+1}$ and hence $D_{a_1(x)}f = 0$ and $D_{a_2(x_s)}g = 0$,
where for $ a \in \mbox{Pin}(n+1)$ the operator $D_a$ is defined in (\ref{Diraca}).

\section{The intertwining formulas for $R_k$ and $R_k^{S}$ and the conformal invariance of $R_k^{S}f=0$}
We can use the intertwining formulas for $D_x$ and $D_{s,x_s}$ given in \cite{LR} to establish the intertwining
formulas for $R_k$ and $R_k^{S}.$

\begin{theorem}\label{ICI}

$$
J_{-1}(C^{-1},x_s)R_{k,u}f(x,u)=R_{k,w}^{S}J(C^{-1},x_s)f(C^{-1}(x_s),\di\frac{(x_s-e_{n+1})w(x_s-e_{n+1})}{\|x_s-e_{n+1}\|^2}),
$$
where $R_{k,u}$ is the Rarita-Schwinger operator in Euclidean space with respect to $u \in \mathbb{R}^{n}$,
$R_{k,w}^{S}$ is the spherical Rarita-Schwinger type operator on $\mathbb{S}^n$ with respect to $w  \in \mathbb{R}^{n+1}$,
$u=\displaystyle\frac{(x_s - e_{n+1})w(x_s - e_{n+1})}{||x_s - e_{n+1}||^2},$
 $J(C^{-1},x_s)=\di\frac{x_s - e_{n+1}}{\|x_s - e_{n+1}\|^n}$
and $J_{-1}(C^{-1},x_s)=\di\frac{x_s - e_{n+1}}{\|x_s - e_{n+1}\|^{n+2}}.$
\end{theorem}
{\bf Proof:}\quad  In \cite{LR} it is shown that $D_x=J_{-1}(C^{-1},x_s)^{-1}D_{s,x_s}J(C^{-1},x_s).$ Consequently,
$$\begin{array}{ll}R_{k,u}f(x,u)=P_{k,u}D_xf(x,u)=P_{k,u}J_{-1}(C^{-1},x_s)^{-1}D_{s,x_s}J(C^{-1},x_s)f(C^{-1}(x_s),u)\end{array}$$
Now applying Theorem \ref{t2}, the previous equation becomes
$$\begin{array}{ll}R_{k,u}f(x,u)=J_{-1}(C^{-1},x_s)^{-1}P_{k,w}D_{s,x_s}J(C^{-1},x_s)f(C^{-1}(x_s),\di\frac{(x_s - e_{n+1})w(x_s - e_{n+1})}{||x_s - e_{n+1}||^2})\\
\\
=J_{-1}(C^{-1},x_s)^{-1}R_{k,w}^SJ(C^{-1},x_s)f(C^{-1}(x_s),\di\frac{(x_s - e_{n+1})w(x_s - e_{n+1})}{||x_s - e_{n+1}||^2}). \qquad \blacksquare
\end{array}
$$

We have the similar result for the Rarita-Schwinger operator under the Cayley transformation.
\begin{theorem}\label{IC}
$$
J_{-1}(C,x)R_{k,u}^Sg(x_s,u)=R_{k,w}J(C,x)g(C(x),\di\frac{(x+e_{n+1})w(x+e_{n+1})}{\|x+e_{n+1}\|^2}),
$$
where $R_{k,u}^S$ is the Rarita-Schwinger type op\-er\-ator on the sphere with re\-spect to $u$
and $R_{k,w}$ is the Rarita-Schwinger op\-er\-a\-tor in Eu\-clidean space with respect to $w$,

\indent$u=\displaystyle\frac{(x + e_{n+1})w(x +e_{n+1})}{||x +e_{n+1}||^2},$
$J(C,x)=\di\frac{x +e_{n+1}}{\|x+e_{n+1}\|^n}$ and $J_{-1}(C,x)=\di\frac{x+e_{n+1}}{\|x +e_{n+1}\|^{n+2}}.$
\end{theorem}

In other words we have the following intertwining relations for $R_k$ and $R_k^S$:
\begin{equation}\label{IR1}
J_{-1}(C^{-1},x_s)R_{k}=R_{k}^{S}J(C^{-1},x_s)
\end{equation}
\begin{equation}\label{IR2}
J_{-1}(C,x)R_{k}^S = R_{k}J(C,x)
\end{equation}

As a corollary to Theorems \ref{ICI} and \ref{IC} we have the conformal invariance of equation $R_{k,w}^S f = 0$:
\begin{theorem}
$ R_{k,u}^Sg(x_s,u) = 0 $ if and only if
$$
R_{k,w} J(C,x)g(C(x),\di\frac{(x+e_{n+1})w(x+e_{n+1})}{\|x+e_{n+1}\|^2})= 0
$$
and $ R_{k,u}f(x,u) = 0 $ if and only if
$$
R_{k,w}^S J(C^{-1},x_s)f(C^{-1}(x_s),\di\frac{(x_s-e_{n+1})w(x_s-e_{n+1})}{\|x_s-e_{n+1}\|^2})= 0.
$$
\end{theorem}

\section{ The fundamental solutions of $R^{S}_k$ and some basic integral formulas}
The reproducing kernel of $\mathcal{M}_k$ with respect to integration over $\mathbb{S}^{n-1}$ is given by
(see \cite{BDS}, \cite{DLRV})
$$
Z_k(u,v):=\di\sum_\sigma P_\sigma(u)V_\sigma(v)v,
$$ where
$$
P_\sigma(u)=\di\frac{1}{k!}\Sigma(u_{i_1}-u_1e_1^{-1}e_{i_1})\ldots(u_{i_k}-u_1e_1^{-1}e_{i_k}),
V_\sigma(v)=\di\frac{\partial^kG(v)}{\partial v_{2}^{j_2}\ldots \partial v_{n}^{j_{n}}}\,
$$
$j_2+\ldots +j_{n}=k, i_k\in\{2,\cdots,n\}, G(v)=\di\frac{-1}{\omega_{n}}\di\frac{v}{\|v\|^{n}},$
and $\omega_{n}$ is the surface area of the unit sphere in $\mathbb{R}^{n}.$ Here summation is
taken over all permutations of the monomials without repetition.  This function is left monogenic
in $u$ and it is a right monogenic polynomial in $v$. It is homogeneous of degree $k$ in both $u$ and $v$. See \cite{BDS} and elsewhere.

Consider the kernel of the Rarita-Schwinger operator in Euclidean $n-$space

$$
\begin{array}{ll}E_k(x-y,u,v)= \di\frac{1}{\omega_{n}c_k}\di\frac{x-y}{\|x-y\|^n}Z_k(\frac{(x-y)u(x-y)}{\|x-y\|^2}, v)\\
\\
=\di\frac{1}{\omega_nc_k}J(C^{-1},x_s)^{-1}\di\frac{x_s-y_s}{\|x_s-y_s\|^n}J(C^{-1},y_s)^{-1}Z_k(\di\frac{(x-y)u(x-y)}{\|x-y\|^2},v),
\end{array}
$$
where $c_k=\di\frac{n-2}{n-2+2k}.$
See for instance \cite{DLRV}. Note that $(x-y)u(x-y)\in \mathbb{R}^{n} $ as $u$, $x $ and $y $ $\in \mathbb{R}^{n}$.

Now applying the Cayley transformation to the above kernel, we obtain

\begin{eqnarray}\label{FunSol1}
&& E_k^S(x_s,y_s,u,v):=
\di\frac{1}{\omega_nc_k}J(C^{-1},x_s)J(C^{-1},x_s)^{-1}\di\frac{x_s-y_s}{\|x_s-y_s\|^n}J(C^{-1},y_s)^{-1}Z_k(au\tilde{a},v)\nonumber\\
&&=\di\frac{1}{\omega_nc_k}\di\frac{x_s-y_s}{\|x_s-y_s\|^n}J(C^{-1},y_s)^{-1}Z_k(au\tilde{a},v),
\end{eqnarray}
where $a=a(x_s,y_s)=\di\frac{J(C^{-1},x_s)^{-1}(x_s-y_s)J(C^{-1},y_s)^{-1}}{\|J(C^{-1},x_s)^{-1}\|\|(x_s-y_s)\|\|J(C^{-1},y_s)^{-1}\|}.$

$E_k^{S}(x_s,y_s, u, v)$ is the fundamental solution to $R_k^{S}f(x_s,u)= 0$ on
$ \mathbb{S}^n $. This function is left monogenic in $u$ and it is also right monogenic in $v.$

In the same way we obtain that
\begin{eqnarray}\label{FunSol2}
\di\frac{1}{\omega_nc_k}Z_k(u,\tilde{a}va)J(C^{-1},y_s)^{-1}\di\frac{x_s-y_s}{\|x_s-y_s\|^n}
\end{eqnarray}
is a non trivial solution to $g(x_s, v)R_{k,r}^{S} = 0$. In fact, this function is $E_k^{S}(x_s,y_s, u, v).$

Applying the same arguments in \cite{DLRV} to prove the representations (\ref{FunSol1}) and (\ref{FunSol2})
are the same up to a reflection, we have

$$
\begin{array}{ll}
\di\frac{1}{\omega_nc_k}Z_k(u,\tilde{a}va)J(C^{-1},y_s)^{-1}\di\frac{x_s-y_s}{\|x_s-y_s\|^n}\\
\\
=-\di\frac{1}{\omega_nc_k}\tilde{a}Z_k(au\tilde{a},v)aJ(C^{-1},y_s)^{-1}\di\frac{x_s-y_s}{\|x_s-y_s\|^n}\\
\\
=-\di\frac{1}{\omega_nc_k}J(C^{-1},y_s)^{-1}\di\frac{x_s-y_s}{\|x_s-y_s\|^n}\di\frac{J(C^{-1},x_s)^{-1}}{\|J(C^{-1},x_s)^{-1}\|}
Z_k(au\tilde{a},v)\di\frac{J(C^{-1},x_s)^{-1}}{\|J(C^{-1},x_s)^{-1}\|}\\
\\
=-\,\di\frac{J(C^{-1},x_s)^{-1}}{\|J(C^{-1},x_s)^{-1}\|}\di\frac{1}{\omega_nc_k}\di\frac{x_s-y_s}{\|x_s-y_s\|^n}J(C^{-1},y_s)^{-1}
Z_k(au\tilde{a},v)\di\frac{J(C^{-1},x_s)^{-1}}{\|J(C^{-1},x_s)^{-1}\|}.
\end{array}
$$
\vspace{0.1cm}
\begin{theorem}\label{stDs} (Stokes' Theorem for the $n$-spherical Dirac operator $D_{s}$) \cite{LR}

Suppose $U_s$ is a domain on $\mathbb{S}^n$ and $f,g: U_s \times \mathbb{R}^{n} \to Cl_{n+1}$
are $C^1$, then for $\partial V_s $ a sufficiently
smooth hypersurface in $U_s$ bounding a subdomain $V_s$ of $U_s$, we have
$$\begin{array}{ll}
\di\int_{\partial V_s} g(x_s, u)n(x_s)f(x_s, u)d\Sigma(x_s)\\
\\
=\di\int_{V_s} (g(x_s, u)D_{s,x_s})f(x_s, u) +  g(x_s, u)(D_{s,x_s}f(x_s, u))dS(x_s) ,
\end{array}$$
where $dS(x_s)$ is the $n$-dimensional area measure on $V_s $, $d\Sigma(x_s)$ is the $n-1$-dimensional scalar
Lebesgue measure on $\partial V_s$ and $n(x_s)$ is the normal vector tangent to the sphere at $x_s$, orthogonal to $\partial V_s$ and
pointing outward.
\end{theorem}
\begin{definition}\label{defiFis} \cite{DLRV} \quad

For any $Cl_{n+1}-$valued polynomials $P(u), Q(u)$,
the inner product $(P(u), Q(u))_{u}$ with respect to $u\in \mathbb{R}^{n}$ is given by
$$
(P(u), Q(u))_{u}=\di\int_{\mathbb{S}^{n-1}} P(u)Q(u)ds(u),
$$
where $\mathbb{S}^{n-1}$ is the unit sphere in $\mathbb{R}^{n}$.
\end{definition}

For any $p_k \in \mathcal{M}_k,$ one obtains
$$
p_k(u)=(Z_k(u,v), p_k(v))_{v} = \int_{\mathbb{S}^{n-1}}Z_k(u,v)p_k(v)ds(v).
$$ See \cite{BDS}.

\begin{theorem}\label{srari}(Stokes' Theorem for the $n$-spherical Rarita-Schwinger type operator $R_k^{S}$)

Let $U_s, V_s, \partial V_s$ be as in Theorem \ref{stDs}. Then for $ f,g \in C^1(U_s \times \mathbb{R}^{n}, {\mathcal M}_{k})$,
we have
\begin{eqnarray*}
\int_{V_s} \Big((g(x_s,u)R_k^{S}, f(x_s,u))_{u} + (g(x_s,u), R_k^{S}f(x_s,u))_{u} \Big)dS(x_s)\\
=\int_{\partial V_s} (g(x_s,u), P_k n(x_s)f(x_s,u))_{u}d\Sigma(x_s),\\
=\int_{\partial V_s} (g(x_s,u)n(x_s)P_{k,r}, f(x_s,u))_{u}d\Sigma(x_s)\\
=\int_{\partial V_s} (g(x_s,u)n(x_s) f(x_s,u))_{u}d\Sigma(x_s)
\end{eqnarray*}
where
$dS(x_s)$ is the $n$-dimensional area measure on $V_s$, $n(x_s)$ and $d\Sigma(x_s)$ are as in Theorem \ref{stDs}.
\end{theorem}
{\bf Proof:}\quad
The proof follows similar lines to the proof of Theorem 6 in \cite{DLRV}.
First, by the traditional Clifford version of Stokes' Theorem 
\begin{eqnarray*}\int_{\partial V_s} (g(x_s,u) n(x_s)f(x_s,u))_{u}d\Sigma(x_s)\\
=\int_{V_s} \Big((g(x_s,u)D_{s,x_s}, f(x_s,u))_{u} + (g(x_s,u), D_{s,x_s}f(x_s,u))_{u} \Big)dS(x_s).\end{eqnarray*}
By applying the Almansi-Fischer decomposition to $g(x_s,u)D_{s,x_s}$ and $D_{s,x_s}f(x_s,u)$ and Definition \ref{defiFis} the right side of the previous equation becomes
$$\int_{V_s} \Big((g(x_s,u)R_k^{S}, f(x_s,u))_{u} + (g(x_s,u), R_k^{S}f(x_s,u))_{u} \Big)dS(x_s).$$
The other identities follow from arguments given in the proof of Theorem 6 in \cite{DLRV}. $\blacksquare$
\begin{coro}(Cauchy's Theorem)

If $R_k^{S}f(x_s,u) = 0$ and $g(x_s, u)R_k^{S} = 0$ for
$f,g\in C^1(U_s \times \mathbb{R}^{n},{\mathcal M}_k)$,
then we have
$$
\int_{\partial V_s} (g(x_s,u), P_kn(x_s)f(x_s,u))_{u}d\Sigma(x_s) = 0,
$$
where $\partial V_s$ is a sufficiently smooth hypersurface in $U_s$ bounding a subdomain $V_s$ of $U_s$.
\end{coro}

Now let us look at Stokes' Theorem for Rarita-Schwinger operators $R_k$ in $\mathbb{R}^n$.\\
Suppose $U$ is a domain on $\mathbb{R}^n$ and $f_{\star},g_{\star}: U \times \mathbb{R}^{n} \to Cl_{n+1}$
are $C^1$, then for $\partial V$ a sufficiently
smooth hypersurface in $U$ bounding a relatively compact subdomain $V$ of $U$, we have
$$\begin{array}{ll}
\di\int_V[(g_{\star}(x,u)R_k, f_{\star}(x,u))_u+(g_{\star}(x,u), R_kf_{\star}(x,u))_u]dx^n\\
\\
=\di\int_{\partial V}\left(g_{\star}(x,u), P_kn(x)f_{\star}(x,u)\right)_ud\sigma(x),
\end{array}$$ where $d\sigma(x)$ is the scalar Lebesgue measure on $\partial V$.
Now consider the integral on the right hand side
$$\begin{array}{lll}\di\int_{\partial V}\left(g_{\star}(x,u),P_kn(x)f_{\star}(x,u)\right)_ud\sigma(x)
=\di\int_{\partial V}\int_{\mathbb{S}^{n-1}}g_{\star}(x,u)P_kn(x)f_{\star}(x,u)ds(u)d\sigma(x)\\
\\
=\di\int_{C(\partial V)}\int_{\mathbb{S}^{n-1}}g_{\star}(C^{-1}(x_s),u)P_{k,u}J(C^{-1},x_s)n(x_s)J(C^{-1},x_s)f_{\star}(C^{-1}(x_s),u)ds(u)d\Sigma(x_s),
\end{array}$$
where $x_s=C(x)$, $C(\partial V)$ bounds a domain $C(V)$ in $\mathbb{S}^n$, $d\Sigma(x_s)$
is the scalar Lebesgue measure on $C(\partial V)$ and $J(C^{-1},x_s)=\di\frac{x_s - e_{n+1}}{\|x_s - e_{n+1}\|^n}.$
Since $P_{k,u}$ in\-ter\-changes with $J(C^{-1},x_s)$, the previous integral becomes
\begin{equation}\label{firsteq}
\begin{array}{ll}
\di\int_{C(\partial V)}\int_{\mathbb{S}^{n-1}}g_{\star}(C^{-1}(x_s),\di\frac{(x_s-e_{n+1})w(x_s-e_{n+1})}{\|x_s-e_{n+1}\|^2})J(C^{-1},x_s)P_{k,w}n(x_s)J(C^{-1},x_s)\\
\\
\qquad \qquad \qquad \qquad \qquad  f_{\star}(C^{-1}(x_s),\di\frac{(x_s-e_{n+1})w(x_s-e_{n+1})}{\|x_s-e_{n+1}\|^2})ds(w)d\Sigma(x_s)
\end{array}
\end{equation}
where $u=\di\frac{(x_s-e_{n+1})w(x_s-e_{n+1})}{\|x_s-e_{n+1}\|^2}.$

Consider the integral on the left hand side
$$\begin{array}{ll}
\di\int_V[(g_{\star}(x,u)R_k, f_{\star}(x,u))_u+(g_{\star}(x,u), R_kf_{\star}(x,u))_u]dx^n\\
\\
=\di\int_V\int_{\mathbb{S}^{n-1}}[g_{\star}(x,u)R_{k,r,u}f_{\star}(x,u)+g_{\star}(x,u)R_{k,u}f_{\star}(x,u)]ds(u)dx^n
\end{array}$$
Applying Theorem 3, the integral now is equal to
$$
\begin{array}{llll}
\di\int_{C(V)}\int_{\mathbb{S}^{n-1}}[g_{\star}(C^{-1}(x_s),\di\frac{(x_s-e_{n+1})w(x_s-e_{n+1})}{\|x_s-e_{n+1}\|^2})
J(C^{-1},x_s)R_{k,r,w}^SJ(C^{-1},x_s)\\
\\
\qquad \qquad \qquad \qquad \qquad \qquad f_{\star}(C^{-1}(x_s),\di\frac{(x_s-e_{n+1})w(x_s-e_{n+1})}{\|x_s-e_{n+1}\|^2})
\end{array}
$$
\begin{equation}\label{secondeq}
\begin{array}{ll}
+g_{\star}(C^{-1}(x_s),\di\frac{(x_s-e_{n+1})w(x_s-e_{n+1})}{\|x_s-e_{n+1}\|^2})J(C^{-1},x_s)R_{k,w}^SJ(C^{-1},x_s)\\
\\
\qquad \qquad \qquad \qquad \qquad f_{\star}(C^{-1}(x_s),\di\frac{(x_s-e_{n+1})w(x_s-e_{n+1})}{\|x_s-e_{n+1}\|^2})]ds(w)dS(x_s)
\end{array}
\end{equation}
where $C(V)=V_s$ is a domain in $\mathbb{S}^n$.

Stokes' Theorem for Rarita-Schwinger operators $R_k$ in $\mathbb{R}^n$ tells us that (\ref{firsteq}) is equal to (\ref{secondeq}).
Therefore Stokes' Theorem for Rarita-Schwinger type operators is conformally invariant under the Cayley transformation.

Now let us consider Stokes' Theorem for $R_k^S$ in $\mathbb{S}^n$.
\begin{eqnarray*}
\di\int_{V_s} \Big((g(x_s,u)R_k^{S}, f(x_s,u))_{u} + (g(x_s,u), R_k^{S}f(x_s,u))_{u} \Big)dS(x_s)\\
=\di\int_{\partial V_s} (g(x_s,u), P_k n(x_s)f(x_s,u))_{u}d\Sigma(x_s),
\end{eqnarray*}
where $V_s, \partial V_s, dS(x_s)$ and $d\Sigma(x_s)$ are stated as in Theorem $7.$

First look at $$\begin{array}{lll}
\di\int_{\partial V_s} (g(x_s,u), P_k n(x_s)f(x_s,u))_{u}d\Sigma(x_s)
=\di\int_{\partial V_s}\int_{\mathbb{S}^{n-1}} g(x_s,u), P_k n(x_s)f(x_s,u)ds(u)\Sigma(x_s)\\
\\
=\di\int_{C^{-1}(\partial V_s)}\int_{\mathbb{S}^{n-1}}g(C(x),u)P_{k,u}J(C(x)n(x)J(C,x)f(C(x),u)ds(u)d\sigma(x),
\end{array}$$
where $J(C,x)=\di\frac{x+e_{n+1}}{\|x+e_{n+1}\|^n}, x=C^{-1}(x_s)$ and $C^{-1}(\partial V_s)$ bounds a domain $C^{-1}(V_s)$ in $\mathbb{R}^n$. Since we can interchange $P_{k,u}$ with $J(C,x),$ the previous integral is equal to

\begin{eqnarray}\label{thirdeq}
\di\int_{C^{-1}(\partial V_s)}\int_{\mathbb{S}^{n-1}}g(C(x),\di\frac{(x+e_{n+1})w(x+e_{n+1})}{\|x+e_{n+1}\|^2})J(C(x)P_{k,w}n(x)J(C,x)\nonumber\\
\nonumber\\
f(C(x),\di\frac{(x+e_{n+1})w(x+e_{n+1})}{\|x+e_{n+1}\|^2})ds(w)d\sigma(x),\end{eqnarray}
where $u=\di\frac{(x+e_{n+1})w(x+e_{n+1})}{\|x+e_{n+1}\|^2}.$

Second we look at
$$\begin{array}{ll}\di\int_{V_s} \Big((g(x_s,u)R_k^{S}, f(x_s,u))_{u} + (g(x_s,u), R_k^{S}f(x_s,u))_{u} \Big)dS(x_s)\\
\\
=\di\int_{V_s}\int_{\mathbb{S}^{n-1}} (g(x_s,u)R_{k,r,u}^{S})f(x_s,u) + g(x_s,u) (R_{k,u}^{S}f(x_s,u))ds(u)dS(x_s).
\end{array}$$
Applying Theorem 4, the integral becomes
\begin{eqnarray}\label{fortheq}
\di\int_{C^{-1}(V_s)}\int_{\mathbb{S}^{n-1}} g(C(x),\di\frac{(x+e_{n+1})w(x+e_{n+1})}{\|x+e_{n+1}\|^2})J(C,x)R_{k,r,w}^{S}\nonumber\\
\nonumber\\
J(C,x)f(C(x),\di\frac{(x+e_{n+1})w(x+e_{n+1})}{\|x+e_{n+1}\|^2})\nonumber\\
 \nonumber\\
 + g(C(x),\di\frac{(x+e_{n+1})w(x+e_{n+1})}{\|x+e_{n+1}\|^2}J(C,x)\nonumber\\
 \nonumber\\
  R_{k,w}^{S}J(C,x)f(C(x),\di\frac{(x+e_{n+1})w(x+e_{n+1})}{\|x+e_{n+1}\|^2})ds(w)dS(x_s).
\end{eqnarray}
Stokes' Theorem for $R_k^S$ on the sphere shows that (\ref{thirdeq}) is equal to (\ref{fortheq}). Thus Stokes' Theorem for Rarita-Schwinger operators is also conformally invariant under the inverse of the Cayley transformation.

\begin{theorem} (Borel-Pompeiu Theorem)
Suppose $U_s$, $V_s$ and $\partial V_s$ are as in Theorem \ref{stDs} and $y_s \in V_s.$
Then for $f \in C^1(U_s \times \mathbb{R}^{n},{\mathcal M}_k) $ we have
\begin{eqnarray*}
f(y_s, u') = J(C^{-1},y_s)\int_{\partial V_s} (E_k^{S}(x_s,y_s, u, v), P_{k}n(x_s) f(x_s,v))_{v}d\Sigma(x_s) \\
-J(C^{-1},y_s)\int_{V_s} (E_k^{S}(x_s,y_s, u, v),R_k^{S}f(x_s,v))_{v} dS(x_s)
\end{eqnarray*}
where $u'=\di\frac{J(C^{-1},y_s)^{-1}uJ(C^{-1},y_s)^{-1}}{\|J(C^{-1},y_s)^{-1}\|^2},$
$dS(x_s)$ is the $n$-dimensional area measure on $V_s \subset \mathbb{S}^n $, $n(x_s)$ and $d\Sigma(x_s)$ are as in Theorem \ref{stDs}.
\end{theorem}

{\bf Proof:}\quad
In this proof we will use the representation of $ E_k^{S}(x_s,y_s, u, v) $  given by (\ref{FunSol2}).

Let $B_s(y_s, \epsilon)$ be the ball centered at $y_s \in \mathbb{S}^n$ with radius $\epsilon$. We denote
$C^{-1}(B_s(y_s, \epsilon))$ by $B(y, r)$, and
$C^{-1}(\partial B_s(y_s, \epsilon))$ by $ \partial B(y, r),$ where $y = C^{-1}(y_s)\in \mathbb{R}^n $ and $r$ is the radius of $B(y, r)$ in $\mathbb{R}^n$. Consider
$ \overline {B}_s(y_s, \epsilon) \subset V_s $, then we have
$$
\begin{array}{ll}
\di\int_{V_s} (E_k^{S}(x_s,y_s, u, v),R_k^{S}f(x_s,v))_{v} dS(x_s) =\\
\\
\di\int_{V_s \setminus B_s(y_s, \epsilon)} (E_k^{S}(x_s,y_s, u, v),R_k^{S}f(x_s,v))_{v} dS(x_s) \\
\\
+\di\int_{B_s(y_s, \epsilon)} (E_k^{S}(x_s,y_s, u, v),R_k^{S}f(x_s,v))_{v} dS(x_s).
\end{array}
$$
Because of the degree of homogeneity of $x_s-y_s$ in $E_k^{S}$, the second integral on the
right hand goes to zero as $\epsilon$ goes to zero.
Applying Theorem \ref{srari} to the first integral on the right hand we obtain
$$
\begin{array}{ll}
\di\int_{V_s \setminus B_s(y_s, \epsilon)}(E_k^{S}(x_s,y_s, u, v),R_k^{S}f(x_s,v))_{v} dS(x_s)\\
\\
=\di\int_{\partial V_s} (E_k^{S}(x_s,y_s, u, v), P_{k}n(x_s) f(x_s,v))_{v}d\Sigma(x_s)\\
\\
-\di\int_{\partial B_s(y_s, \epsilon)}(E_k^{S}(x_s,y_s, u, v), P_{k}n(x_s) f(x_s,v))_{v}d\Sigma(x_s).
\end{array}
$$
Since $f(x_s, v) = (f(x_s, v) - f(y_s, v)) + f(y_s, v)$ and taking into account the degree of
homogeneity of $x_s-y_s$ in $E_k^{S}$ and the continuity  of $f$, we can replace the second integral on the right hand by
$$
\di\int_{\partial B_s(y_s, \epsilon)} (E_k^{S}(x_s,y_s, u, v), P_{k} n(x_s) f(y_s,v))_{v}d\Sigma(x_s).
$$
Applying Theorem \ref{srari}, this integral is equal to
$$
\begin{array}{lll}
\di\int_{\partial B_s(y_s, \epsilon)} (E_k^{S}(x_s,y_s, u, v), n(x_s) f(y_s,v))_{v}d\Sigma(x_s)\\
\\
=\di\int_{\partial B_s(y_s, \epsilon)} \int_{\mathbb{S}^{n-1}}E_k^{S}(x_s,y_s, u, v) n(x_s) f(y_s,v)d\Sigma(x_s) ds(v)\\
\\
=\di\int_{\partial B_s(y_s, \epsilon)} \int_{\mathbb{S}^{n-1}}\di\frac{1}{\omega_nc_k}
Z_k(u,\tilde{a}va)
J(C^{-1},y_s)^{-1}\di\frac{x_s-y_s}{\|x_s-y_s\|^n}n(x_s) f(y_s,v)ds(v)d\Sigma(x_s).
\end{array}
$$
Now applying the inverse of the Cayley transformation to the last integral, we have
$$
\begin{array}{lll}
\di\int_{\partial B(y, r)} \int_{\mathbb{S}^{n-1}}
\di\frac{1}{\omega_nc_k}Z_k(u,\di\frac{(x-y)w(x-y)}{\|x-y\|^2})
J(C^{-1},y_s)^{-1}J(C,y)^{-1}\di\frac{x-y}{\|x-y\|^n}\\
\\
J(C,x)^{-1}J(C,x)n(x)J(C,x)f(C(y),\di\frac{J(C,y)wJ(C,y)}{\|J(C,y)\|^2})ds(w)d\sigma(x),
\end{array}
$$
where $d\sigma(x)$ is the $n-1-$ dimensional scalar Lebesgue measure on $\partial B(y,r)$ in $\mathbb{R}^n$ and
$v=\di\frac{J(C,y)wJ(C,y)}{\|J(C,y)\|^2}, w\in \mathbb{R}^{n+1}.$
In fact, $v$ is a vector in $\mathbb{R}^{n}$ which is obtained by reflecting $w$ in $\mathbb{R}^{n+1}$ and its last component is a constant.

Place $J(C,x)=(J(C,x)-J(C,y))+J(C,y)$, but $J(C,x)-J(C,y)$ tends to zero as $x$ approaches $y$.
Thus the previous integral can be replaced by
$$
\begin{array}{lll}
\di\int_{\partial B(y, r)} \int_{\mathbb{S}^{n-1}}
\di\frac{1}{\omega_nc_k}Z_k(u,\di\frac{(x-y)w(x-y)}{\|x-y\|^2})
\di\frac{x-y}{\|x-y\|^n}n(x)\\
\\
J(C,y)f(C(y),\di\frac{J(C,y)wJ(C,y)}{\|J(C,y)\|^2})ds(w)d\sigma(x).
\end{array}
$$

Here $n(x)=\di\frac{y-x}{\|x-y\|}$ is the unit out normal vector. Now the last integral becomes
$$\begin{array}{ll}
\di\int_{\partial B(y, r)} \di\int_{\mathbb{S}^{n-1}}
\di\frac{1}{\omega_nc_k}Z_k(u,\di\frac{(x-y)w(x-y)}{\|x-y\|^2})\di\frac{1}{\|x-y\|^{n-1}}\\
\\
J(C,y)f(C(y),\di\frac{J(C,y)wJ(C,y)}{\|J(C,y)\|^2})ds(w)d\sigma(x),
\end{array}$$
Using Lemma $5$ in \cite{DLRV}, the integral is now
$$
\begin{array}{lll}
\di\int_{\mathbb{S}^{n-1}} Z_k(u, w)J(C,y)f(C(y),\di\frac{J(C,y)wJ(C,y)}{\|J(C,y)\|^2})ds(w) \\
\\
= J(C,y)f(C(y),\di\frac{J(C,y)uJ(C,y)}{\|J(C,y)\|^2})=J(C^{-1},y_s)^{-1}f(y_s, \di\frac{J(C^{-1},y_s)^{-1}uJ(C^{-1},y_s)^{-1}}{\|J(C^{-1},y_s)^{-1}\|^2}),
\end{array}
$$
since $J(C,y)=J(C^{-1},y_s)^{-1}.$

Now by setting $u'=\di\frac{J(C^{-1},y_s)^{-1}uJ(C^{-1},y_s)^{-1}}{\|J(C^{-1},y_s)^{-1}\|^2}$ and multiplying both sides of the above equation by $J(C^{-1},y_s),$ we obtain

$$f(y_s,u')=J(C^{-1},y_s)\di\int_{\mathbb{S}^{n-1}} Z_k(u, w)J(C,y)f(C(y),\di\frac{J(C,y)wJ(C,y)}{\|J(C,y)\|^2})ds(w).$$

Therefore when $\epsilon$ tends to zero we get the desired result. ~ $\blacksquare$

\begin{coro}

Let $\Psi$ be a function in $C^\infty (V_s \times \mathbb{R}^{n},\mathcal{M}_k)$
and $\mbox{supp}(\Psi)\subset V_s $. Then
$$
\Psi(y_s, u') = - J(C^{-1}, y_s)\int_{V_s} (E_k^{S}(x_s,y_s, u, v),R_k^{S}\Psi(x_s,v))_{v} dS(x_s),
$$
where $u'=\di\frac{J(C^{-1},y_s)^{-1}uJ(C^{-1},y_s)^{-1}}{\|J(C^{-1},y_s)^{-1}\|^2}$.
\end{coro}

\begin{coro} (Cauchy Integral Formula for $R_k^{S}$)

If $R_k^{S}f(x_s, u) = 0$, then for $y_s \in V_s$ we have
\begin{eqnarray*}
f(y_s, u') &=&J(C^{-1}, y_s)\int_{\partial V_s} (E_k^{S}(x_s,y_s, u, v), P_kn(x_s)f(x_s,v))_{v} d\Sigma(x_s) \\
&=&J(C^{-1}, y_s)\int_{\partial V_s} (E_k^{S}(x_s,y_s, u, v) n(x_s)P_{k,r},f(x_s,v))_{v} d\Sigma(x_s),
\end{eqnarray*}
where $u'=\di\frac{J(C^{-1},y_s)^{-1}uJ(C^{-1},y_s)^{-1}}{\|J(C^{-1},y_s)^{-1}\|^2}$.
\end{coro}

\begin{definition} (Cauchy Transform for $R_k^{S}$)

For a domain $V_s \subset \mathbb{S}^n $ and a function $f(x_s,u): V_s \times \mathbb{R}^{n} \to Cl_{n+1},$ which is monogenic in $u,$
the $T_k$-transform of $f$ is defined to be
$$
(T_kf)(y_s, v) = - \int_{V_s} (E_k^{S}(x_s,y_s, u, v), f(x_s, u))_{u} dS(x_s), ~~ \mbox{for}~~ y_s \in V_s.
$$
\end{definition}

\begin{theorem}

For a function $\psi$ in $C^\infty (\mathbb{S}^n \times \mathbb{R}^{n},\mathcal{M}_k)$ we have
$$
P_k J(C^{-1},y_s)D_{s,y_s}\int_{\mathbb{S}^{n}} -(E_k^{S}(x_s,y_s, u, v), \psi(x_s, u))_{u} dS(x_s) = \psi(y_s, v).
$$
\end{theorem}
{\bf Proof:} \quad By \cite{LR}, the integral
$$P_k J(C^{-1},y_s)D_{s,y_s}\int_{\mathbb{S}^{n}} -(E_k^{S}(x_s,y_s, u, v), \psi(x_s, u))_{u} dS(x_s)$$
can be replaced by
$$
P_k J(C^{-1},y_s)\int_{\partial B_s(y_s,\epsilon)} -n(x_s)(E_k^{S}(x_s,y_s, u, v), \psi(x_s, u))_{u} dS(x_s),$$
which in turn is equal to
$$\begin{array}{lll}P_k J(C^{-1},y_s)\di\int_{\partial B_s(y_s,\epsilon)}\int_{\mathbb{S}^{n-1}} -n(x_s)\di\frac{1}{\omega_nc_k }\di\frac{x_s-y_s}{\|x_s-y_s\|^n}J(C^{-1},y_s)^{-1}Z_k(au\tilde{a},v) \psi(x_s, u)ds(u) dS(x_s).\end{array}$$
Since $\psi(x_s, u)=\psi(x_s, u)-\psi(y_s, u)+\psi(y_s, u)$ then using the continuity of $\psi$, we can replace the previous integral by

$$\begin{array}{lll}P_k J(C^{-1},y_s)\di\int_{\partial B_s(y_s,\epsilon)}\int_{\mathbb{S}^{n-1}} -n(x_s)\di\frac{1}{\omega_nc_k }\di\frac{x_s-y_s}{\|x_s-y_s\|^n}J(C^{-1},y_s)^{-1}Z_k(au\tilde{a},v) \psi(y_s, u)ds(u) dS(x_s).\end{array}$$
Now applying the inverse of the Cayley transformation to the previous integral it becomes
$$\begin{array}{ll}
P_k J(C^{-1},y_s)\di\int_{\partial B(y,r)}\int_{\mathbb{S}^{n-1}} -J(C,x)n(x)J(C,x)\di\frac{1}{\omega_nc_k }J(C,x)^{-1}\di\frac{x-y}{\|x-y\|^n}J(C,y)^{-1}J(C^{-1},y_s)^{-1}\\
\\
Z_k(\di\frac{(x-y)u(x-y)}{\|x-y\|^2},v) \psi(C(y), \di\frac{J(C,y)wJ(C,y)}{\|J(C,y)\|^2})ds(u) d\sigma(x)\\
\\
=P_k J(C^{-1},y_s)\di\int_{\partial B(y,r)}\int_{\mathbb{S}^{n-1}} -J(C,x)n(x)\di\frac{1}{\omega_nc_k }\di\frac{x-y}{\|x-y\|^n}Z_k(\di\frac{(x-y)u(x-y)}{\|x-y\|^2},v)\\
\\
 \psi(C(y), \di\frac{J(C,y)wJ(C,y)}{\|J(C,y)\|^2})ds(u) d\sigma(x),
\end{array}$$ where $u=\di\frac{J(C,y)wJ(C,y)}{\|J(C,y)\|^2}.$

Using the fact $J(C,x)=(J(C,x)-J(C,y))+J(C,y)$, and $J(C,x)-J(C,y)$ tends to zero as $x$ approaches $y,$
the integral can be replaced by
$$\begin{array}{llll}
P_k J(C^{-1},y_s)\di\int_{\partial B(y,r)}\int_{\mathbb{S}^{n-1}} -\di\frac{1}{\omega_nc_k }J(C,y)n(x)\di\frac{x-y}{\|x-y\|^n}Z_k(\di\frac{(x-y)u(x-y)}{\|x-y\|^2},v)\\
\\
 \psi(C(y), \di\frac{J(C,y)wJ(C,y)}{\|J(C,y)\|^2})ds(u) d\sigma(x)\\
\\
=P_k \di\int_{\partial B(y,r)}\int_{\mathbb{S}^{n-1}} \di\frac{1}{\omega_nc_kr^{n-1} }
Z_k(\di\frac{(x-y)u(x-y)}{\|x-y\|^2},v) \psi(C(y), \di\frac{J(C,y)wJ(C,y)}{\|J(C,y)\|^2})ds(u) d\sigma(x)
\end{array}
$$
Applying Lemma 5 in \cite{DLRV}, the integral becomes
$$\begin{array}{llll}
=P_k \di\int_{\mathbb{S}^{n-1}}Z_k(u,v) \psi(C(y), \di\frac{J(C,y)wJ(C,y)}{\|J(C,y)\|^2})ds(u) \\
\\
=P_k \di\int_{\mathbb{S}^{n-1}}Z_k(u,v) \psi(C(y), u)ds(u)=P_k\psi(C(y),v)=\psi(y_s, v). ~\blacksquare
\end{array}
$$

\section{Rarita-Schwinger type operators on real projective space}

We consider $\mathbb{S}^n$ and $\Gamma=\{\pm 1\}$, then $\mathbb{S}^n/\Gamma$ is $\mathbb{R}P^n$, the real projective space. In all that follows $\mathbb{S}^n$ will be a universal covering space of the conformally flat manifold $\mathbb{R}P^n.$ So there is a projection map $p:\mathbb{S}^n \to \mathbb{R}P^n.$ Further for each $x\in \mathbb{S}^n$ we shall denote $p(x)$ by $x'.$ Furthermore if $Q$ is a subset of $U$ then we denote $p(Q)$ by $Q'.$

Consider the trivial bundle $\mathbb{S}^n\times Cl_{n+1}$, then we set up a spinor bundle $E_1$ over $\mathbb{R}P^n$ by making the identification of $(x,X)$ with $(-x,X)$, where $x\in \mathbb{S}^n$ and $X\in Cl_{n+1}$.

Now we change the spherical Cauchy kernel $G_{S}(x-y)=\di\frac{-1}{\omega_n}\di\frac{x-y}{\|x-y\|^n}, x,y \in\mathbb{S}^n,$ for the spherical Dirac operator into a
kernel which is invariant with respect to $\{\pm 1\}$ in the variable $x\in \mathbb{S}^n$. Hence we consider $G_{S}(x-y)+G_{S}(-x-y)$. See \cite{KR}.

Suppose  $V$ is a domain lying in the open northern hemisphere.
We assume
$$
f(x,u): V\times \mathbb{R}^{n} \to Cl_{n+1}
$$ is a $C^1$ function in $x$ and monogenic in $u$.
We observe that the projection map $p: \mathbb{S}^n \to \mathbb{R}P^n$ induces a well defined function
$$
f'(x', u): V' \times \mathbb{R}^{n} \to E_1
$$
such that $f'(x', u) = f(p^{-1} (x'), u)$, where $V'$ is a well defined domain in $\mathbb{R}P^n$ and $x' = p(x)$.

We define the Rarita-Schwinger type operators on $\mathbb{R}P^n$, which we will call the real projective Rarita-Schwinger type operators, in the following form
$$
R_k^{RP^n}f'(x',u)=P_k D_{RP^n,x'} f'(x',u),
$$
where $D_{RP^n,x'}$ is the Dirac operator in $ \mathbb{R}P^n $ with respect to the variable $x'$. See \cite{KR}.

Now we introduce the spherical Rarita-Schwinger kernel which is also invariant with respect to $\{\pm 1\}$ in the variable $x\in \mathbb{S}^n:$
$$
\mathcal{E}_k^{S,1}(x,y,u,v):=E_k^{S}(x,y,u,v)+E_k^{S}(-x,y,u,v).
$$

Through the projection map $p$ (over $x,y\in \mathbb{S}^n$) we obtain a kernel $\mathcal{E}_k^{RP^n,1}(x',y',u,v)$ for $\mathbb{R}P^n$ defined by
$$
\mathcal{E}_k^{RP^n,1}(x',y',u,v)=\mathcal{E}_k^{S,1}(p^{-1}(x'),p^{-1}(y'),u,v).
$$

Now suppose that $S$ is a suitably smooth hypersurface lying in  the open northern hemisphere of
$\mathbb{S}^n$ bounding a subdomain $W$ of $V$ with closure of $W$ in $V.$
\begin{theorem} \label{thm10} If $R_k^{S}f(x,u) = 0$ then for $y\in W$
$$
f(y,w) = J(C^{-1},y)\int_S ( \mathcal{E}_k^{S,1}(x,y,u,v), P_k n(x)f(x,u))_u d\Sigma(x),
$$
where $w=\di\frac{J(C^{-1},y)^{-1}vJ(C^{-1},y)^{-1}}{\|J(C^{-1},y)^{-1}\|^2},$ $n(x)$ is the unit outer normal vector to $S$ at $x$ lying in the tangent space
of $S^n$ at $x$ and $\Sigma$ is the usual Lebesgue measure on $S$.
\end{theorem}
Due to the projection map we have also
\begin{theorem}\label{thm11}
$$
f'(y',\hat{v}) =  J(C^{-1},y')\int_{S'} ( \mathcal{E}_k^{RP^n,1}(x',y',u,v), P_k dp(n(x))f'(x',u))_u d\Sigma'(x'),
$$
where $\hat{v}=\di\frac{J(C^{-1},y')^{-1}vJ(C^{-1},y')^{-1}}{\|J(C^{-1},y')^{-1}\|^2}$, $x' = p(x)$, $y'= p(y)$ and $S'$ is the projection of $S$. Further $\Sigma'$ is a induced
measure on $S'$ from the measure $\Sigma$ on $S$ and $dp$ is the derivative of $p$.
\end{theorem}

Now we will assume that the domain $V$ is such that $-x \in V$ for each $x \in V$ and the
function $f$ is two fold periodic, so that $f(x) = f(-x)$.
Now the projection map $p$ gives rise to a well defined domain $V'$ on $\mathbb{R}P^n$ and a
well defined function $f'(x',u): V'\to E_1$ such that $f'(x',u)=f(\pm x,u)$ for
$p(\pm x)=x'.$ Then if $R_k^{RP^n}f'(x',u) = 0$,
we also have
$$
f'(y',\hat{v}) =  J(C^{-1},y')\int_{S'} ( \mathcal{E}_k^{RP^n,1}(x',y',u,v), P_k dp(n(x))f'(x',u))_u d\Sigma'(x'),
$$ where $\hat{v}$ is stated as in Theorem \ref{thm11}.

If now we suppose that the hypersurface $S$ satisfies $-S = S$ then
both $y$ and $-y$ belong to the subdomain $V$ and in this case
$$
 J(C^{-1},y')\int_{S'} ( \mathcal{E}_k^{RP^n,1}(x',y',u,v), P_k dp(n(x))f'(x',u))_u d\Sigma'(x')= 2 f'(y',\hat{v}).
$$

We can also construct a second spinor bundle $E_2$ over $\mathbb{R}P^n$ by making the identification of
$(x,X)$ with $(-x,-X)$, where $x\in \mathbb{S}^n$ and $X\in Cl_{n+1}$, we introduce the kernel:
$$
\mathcal{E}_k^{S,2}(x,y,u,v):=E_k^{S}(x,y,u,v)-E_k^{S}(-x,y,u,v).
$$
This kernel induces through the projection map on the variable $x,y\in \mathbb{S}^n,$ the kernel on $\mathbb{R}P^n$
$$
\mathcal{E}_k^{RP^n,2}(x',y',u,v)=\mathcal{E}_k^{S,2}(p^{-1}(x'),p^{-1}(y'),u,v).
$$
In this case a solution of Rarita-Schwinger type equation on $\mathbb{R}P^n$
$$
f'(x', u): V' \times \mathbb{R}^{n} \to E_2
$$
will lift to a solution of spherical-Rarita-Schwinger type equation:
 $ f(x, u): V \times \mathbb{R}^{n} \to Cl_{n+1}  $
such that $f(x, u) = - f(-x, u)$.

Suppose that $V$ as before is a domain on $\mathbb{S}^n$ and $S$ is a hypersurface  in $V$ bounding
a subdomain $W$ of $V$. Suppose further that $ f(x, u): V \times \mathbb{R}^{n} \to Cl_{n+1}$ is a solution of the
spherical Rarita-Schwinger type equation such that $f(x, u) = - f(-x, u)$. If $S$ lies entirely in one open hemisphere then
$$
f(y,w) = J(C^{-1},y) \int_S ( \mathcal{E}_k^{S,2}(x,y,u,v), P_k n(x)f(x,u))_u d\Sigma(x),
$$
for each $y \in W$, where $w=\di\frac{J(C^{-1},y)^{-1}vJ(C^{-1},y)^{-1}}{\|J(C^{-1},y)^{-1}\|^2}.$

Via the projection $p$ this integral formula induces the following
$$
f'(y',\hat{v}) =  J(C^{-1},y')\int_{S'} ( \mathcal{E}_k^{RP^n,2}(x',y',u,v), P_k dp(n(x))f'(x',u))_u d\Sigma'(x'),
$$where $\hat{v}$ is stated as in Theorem \ref{thm11}.

On the other hand if $S$ is such that $S = -S$ then
$$
\int_S ( \mathcal{E}_k^{S,2}(x,y,u,v), P_k n(x)f(x,u))_u d\Sigma(x) = 0.
$$

Junxia Li \quad Email:  jxl004@uark.edu\\
John Ryan \quad Email: jryan@uark.edu\\
Carmen J. Vanegas \quad Email: cvanegas@usb.ve

\end{document}